\documentclass[12 pt]{article}

\usepackage{inputenc}
\usepackage{amssymb,amsmath,latexsym}
\usepackage{bm}
\usepackage{graphicx}

\newcommand{\QED}{~\rule{7pt}{7pt}}

\begin{document}

{\Large
\centerline{The polynomial algorithm for optimal spanning hyperforest problem}
\vspace{3mm}
\centerline{Abakarov A.Sh.\footnote{alik@vega.math.spbu.ru}, Sushkov Yu. A.\footnote{Yu.Sushkov@pobox.spbu.ru}}
}

\vspace{5mm}

\centerline{\bf Abstract}

This paper is devoted to one theory of hypergraph connectivity and presents the proof of the polynomial 
algorithm for finding optimal spanning hyperforest~(hypertree) for any given weighed $q$-uniform
hypergraph . 

\vspace{5mm}

{\bf \large Introduction}

\vspace{2mm} Couple $\Gamma = \langle Z, D \rangle$ is a hypergraph, where $Z$ is finite set of 
vertexex, and $D$  is set of edges and any edge $a \in D$ is  subset of $Z$.

Except graph notation, the symbol $\Gamma$ is additionally used to define the following mapping:
$$
\Gamma : D \rightarrow 2^Z,              
$$
where $2^Z$ is a family of all possible sets of $Z$~[3,4,5].

In particular, if suppose that
$$
\Gamma : D \rightarrow {Z \choose q} \subset 2^Z,      \eqno(1)
$$
where ${Z \choose q}$ is a family of all possible $q$-elements subset of $Z$, then the 
hypergraph $\Gamma = \langle Z, D \rangle$ corresponding the mapping~(1) is called $q$-uniform hypergraph.

Hereinafter, only $q$-uniform hypergraphs are considered  and thus, instead this term either  
terms $q$-graph or graph is used.

{\bf Definition 1.} {\it Hyperforest} $\langle \Gamma W, W \rangle$, is  hypergraph where 
for any $A \subseteq W$ the following {\it condition of independence}
$$
|\Gamma A| - q + 1 \ge |A|,                   \eqno(2)
$$
is satisfyed.

Moreover, if the following equation
$$
|\Gamma W| - q + 1 = |W|,                     \eqno(3)
$$
is satisfied for hyperforest $\langle \Gamma W, W \rangle$ such kind of hyperforest 
is referred to as {\it hypertree}.

Set of edges satisfies the condition~(2) is reffered to as {\it independent one}. Therefore, 
any hypergraph edges subset of $D$ spanning the hyperforest is independent.

{\bf Definition 2.} Subgraph  $\langle \Gamma W, W \rangle$ of graph 
$\Gamma = \langle \Gamma D, D \rangle$ with $ \Gamma W = \Gamma D$ and with 
maximum $W$ is referred to as {\it skeleton} of hypergraph $\langle \Gamma D, D \rangle$.

In general case hypergraph skeleton  is hyperforest.

Work~[4] proved that for the skeleton $\langle \Gamma W, W \rangle$ of graph 
$\langle \Gamma D, D \rangle$ exist only one subgraph decomposition $\langle \Gamma T_i, T_i \rangle$, 
such that:

 1) $W = T_1 + T_2 + \ldots + T_b$, $\Gamma T_i = |T_i| + q - 1$, $\forall i\in 1:b$,

 2) for any such $A \in W$, that $\Gamma A = |A| + q - 1 $, there exist only one $i \in 1:b$, that $A \in T_i$.

{\bf Definition 3.} Subgraphs  $\langle \Gamma T_i, T_i \rangle$, $i \in 1:b$, spanned by $W$ set $T_i$ 
partition elements  are referred to as {\it connection components} of  hyperforest $\langle \Gamma W, W \rangle$.

It is obviously that a hyperforest $\langle \Gamma W, W \rangle$ consisted  of only one connected component 
is a hypertree.

If hypergraph $D$ set  $\Gamma =  \langle \Gamma D, D \rangle$ is dependable one the hypergraph can 
have more than one skeleton.
Connection components of any skeletons of this graph univalently induces the connection 
components $ \langle \Gamma D_i, D_i \rangle$, $T_i \subseteq D_i$, $i\in 1:b$, of any hypergraph
$\Gamma $, without any relation to selected skeleton~[4].

Let   $W = T_1 + T_2 + \ldots + T_b$ is graph edges set partition to connection components 
of any hypegraph skeleton $\Gamma = \langle \Gamma D, D \rangle$, $W \subseteq D$, 
while \mbox {$D = D_1 + D_2 + \ldots + D_b$} is partition of its edges, where
 $T_i \subseteq D_i$, $i\in 1:b$. Thus, an edge $a \in D_i \setminus T_i$  is referred
to as {\it graph link} of connection component $\langle \Gamma T_i, T_i \rangle$ with respect 
to the specified graph skeleton.

Therefore, if $H_i$ is connection components graph links set of $T_i$, i.e. $D_i = T_i + H_i$, 
$i\in 1:b$, then hypergraph $\Gamma = \langle \Gamma D, D \rangle$ links set partition $D = D_1 + D_2
+ \ldots + D_b$ does not depend on graph skeleton choice.

{\bf Definition 4.} Subgraph  $\langle \Gamma D_i, D_i \rangle$, $i \in 1:b$, provided by 
partition element $D_i$ of $D$ set is referred to as hypergraph $\langle \Gamma D, D \rangle$ {\it
connection component}.

Let us consider the aforementioned definitions using the $3$-graph represented in the Figure.~1~ (this hypergraph is designated as $\tilde \Gamma = \langle \tilde \Gamma D, D \rangle$).


Hypergraph $\tilde \Gamma$ includes two skeletons: $\langle \tilde \Gamma W_1, W_1 \rangle$, 
$W_1 = \{a,b,c,d,e,g,h\}$ and $\langle \tilde \Gamma W_2, W_2 \rangle$, $W_2 = \{a,b,c,d,f,g,h\}$ which
are hypertrees, i.e. hypergraph $\tilde \Gamma$ edges set $D$ forms one connection component 
$\langle \tilde \Gamma D_1, D_1 \rangle$, $D_1 = D$. Edge $f$ is connection component  
$\langle \tilde \Gamma T_1, T_1 \rangle$, $T_1 = W_1$ link, while edge $e$ is link 
of  $\langle \tilde \Gamma W_2, W_2 \rangle$.

If remove the edges $e$ and $f$ from rib set $D$ of hypergraph $\tilde \Gamma$ then 
the $\tilde \Gamma$ consists only the spanning set $\langle \tilde \Gamma W, W \rangle$, 
$W = \{a,b,c,d,g,h\}$, which is hyperforest uniquely separated in two connection 
components : $\langle \tilde \Gamma T_1, T_1 \rangle$, $T_1 = \{a,b,c\}$ 
and $\langle \tilde \Gamma T_2, T_2 \rangle$, $T_2 = \{d,g,f\}$.

\vspace{5mm}
{\bf \large Optimal skeleton development problem}

Let real function
$$ \omega : D \rightarrow \mathbb{R}^+,
\eqno(4)
$$
is defined in the set of edges of hypergraph  $\Gamma = \langle Z, D \rangle$ 
where the value of $\omega (d)$, $d \in D$ is referred to as  {\it weight} of element $d$. 
In this case, optimal  skeleton development problem for hypergraph  $\Gamma = \langle Z, D \rangle$
is to determine the skeleton $\langle \Gamma W, W \rangle$ with the minimal value of function
$$
\omega(W) = \sum_{d \in W} \omega(d).
$$

Articles~[4,5] show that optimal skeleton development problem is a particular case of 
minimal (maximal) matroid independent set definition problem.

Let us introduce necessary definitions ~[2].

{\bf Definition 5.} Let $D$ is nonvacuous finite set and $I \subseteq 2^D$ is nonvacuous 
subset family of $D$ which meets the following requirements:

1) if $A \subseteq B \in I$ then $A \in I$;

2) if $A \in I$, $B \in I$ and $|A| = |B| - 1$ then exists an
element $a\in B \setminus A$ such as $A + a \in I$.

Then couple $M = \langle I, D \rangle$ is referred to as {\it
matroid}.

Subset family $I$ elements are referred to as {\it independent} set and elements of set 
$2^D \setminus I$ are referred to as  {\it dependent} sets of matroid $M$. Independent sets of family ~$I$
with ultimate number of elements are referred to as {\it basis} of matroid ~$M$.

Matroid can also be defined with the help of  $ \varphi : 2^D \to \mathbb{Z} $, $ \varphi (\emptyset) = 0$ 
function which takes the whole-number values on the subsets of $D$ and satisfies the
following conditions:

1) if  $A \subseteq B \subseteq D$, then
$$
 \varphi(A) \le \varphi(B),  
$$

2) for any two subsets $A, B  \subseteq D$
$$
\varphi (A) + \varphi (B) \ge  \varphi (A \cup B) + \varphi (A \cap B).
$$

While taking the aforementioned in consideration, the set $A$ of matroid 
$M = \langle I, D \rangle$ is  {\it independent } one ($A \in I$) if and only if for any
любого $ A \subseteq D $
$$
 \varphi (A) \ge |A|.             \eqno (5)
$$

This equation is referred to as  {\it independence} condition.

If function included in  equation (5) is taken as
$$
 \varphi (A) = |\Gamma A| - q + 1,
$$
then inequations (2) and (5) are agreed. Hence, matroid $M =
\langle I, D \rangle$ can be connected with any hypergraph
 $\Gamma = \langle Z, D
\rangle$ in case all possible independent hypergraph ~$\Gamma$ edges $D$ subsets are 
taken as independent subsets~$I$ of matroid $M$. The matroid obtained by aforemntioned way is referred to as
 {\it graph matroid}. All possible hypergraph $\Gamma$ skeletons are referred to as  graph matroid {\it
bases} .

It is known~[2] that to allocate the independent matroid subset with the minimal weight, the so called {\it greedy} algorithm can be used. This algorithm can be described as follows:

\vspace{3mm} \underline{\bf Greedy algorithm.} \vspace{2mm}

{\bf Begin}

{\bf Step 1.} Arrange the set $D$ ascending by weights in such a
way as

\hspace*{15mm} $D = \{ d_1, d_2, \ldots, d_{|D|} \}, \ \ \mbox{где} \ \ \omega(d_1) \le
\omega(d_2) \le \omega(d_{|D|})$.

{\bf Step 2.} $W := \emptyset$.

{\bf Step 3.} For any $i \in 1:|D|$ do the following:

\hspace*{15mm} If $W \cup \{d_i\} \in I$, then $W := W \cup
\{d_i\}$;

\hspace*{15mm} else  $W := W$.

{\bf End}

Having the greedy algorithm work completed the independent set  $W \subseteq D$ is built. This set is the matroid
$M$ basis provided with the minimal weight $\omega(W)$.

The described greedy algorithm complexity estimation depends on step~3 efficiency. This step is  running subset
 $W \cup \{d_i\} = W +d_i$ checking for independent subset family  $I$ belonging.

In case of graph matroid, this step checks the correctness of independence condition~(5) 
for any subset $A \subseteq W$ of the hypergraph. It can be done using the complete enumeration of all
possible edges $A \subseteq W$ subset.

In case of graph matroid it can be found that the greedy algorithm step~3 can be performed for polynomial time. Therefore, optimal hypergraph skeleton development algorithm is the polynomial one,
i.e. optimal skeleton development problem belongs to polynomial problems class.

It is obvious that optimal skeleton development problem of $2$-graphs belongs the same class. Prim and Kruskal~[2] algorithms can be used to solve this problem.

It is known that for any hypergraph $\Gamma = \langle Z, D \rangle$ its 
{\it Koenig representation} $\langle D, Z, \Gamma \rangle$ can be assigned. This representation 
is bichromatic graph with two sets of vertexes  $D$, $Z$ and with set of edges $\{(a,b)
\ | \ a\in D \ \ \& \ \ b\in \Gamma a \subseteq Z \}$.

Hypergraph skeleton development problem is closely connected with its Koenig 
representation  {\it complete matching} development, i.e. development of graph  
$\langle D, Z, \Gamma \rangle$ subgraph with the maximal number of edges and with degree of all vertexes
equal to identity.

Hypergraph Koeing representation  $\tilde \Gamma$ is shown in the Figure~2.



{\bf Theorem 1.} {\it Let $\langle  D, Z, \Gamma \rangle$ is the Koenig representation 
of   $q$-uniform hypergraph $\Gamma = \langle Z, D\rangle$. Then removing  any vertexes set $B$ such as
\ $|B| = q-1$, from set $Z$ of graph  $\langle D, Z, \Gamma \rangle$ in the obtained bichromatic graph  $\langle D, Z
\setminus B, \Gamma \rangle$   the complete matching exists if and only if the hypergraph
$\Gamma = \langle Z, D\rangle$ is the hyperforest  (edges set $D$ is independent one)}.

{\bf Proof.}
 The aforementioned means that for any subset $ A
\subseteq D$ the following is true:
$$
|\Gamma A| - (q - 1) \ge |A|.      \eqno(6)
$$

Let us firstly proof the "only if" case, i.e. if for the hypergraph $\Gamma = \langle Z, D\rangle$ 
condition~(6) is fulfilled then complete matching exists in its Koenig representation while 
removing any two   $q-1$ vertexes of $Z$.

For any $A \subseteq D$ the inequality $|\Gamma A| \ge |A|$ is true no matter what 
vertexes $q-1$ belonging the $Z$ are removed. This statement is true due to Hall theorem.

Let us proof the "if" case: if remove any $|B| = q - 1$ vertexes of graph  $\langle  D, \ Z\setminus B,\  \Gamma \rangle$ from the set $Z$  the   Koenig representation complete matching exists and the hypergraph  $\Gamma = \langle Z, D\rangle$ is the hyperforest (edges set $D$ meets the condition (6)).

In order to prove this conclusion let us suppose an opposite, i.e. while removing the set $B$ from hypergraph Koenig representation $\Gamma = \langle Z, D\rangle$ vertexes set $Z$ the obtained graph
$\langle D, \ Z\setminus B, \ \Gamma \rangle$ does not include the complete matching. That means exists the such $A \subseteq D$ as $|\Gamma A| < |A|$. However this means that in case return to hypergraf Koenig representation of vertexes set  $B$ of subset $A \subseteq D$ the independence condition~(6)does not meet.

Thus, the contradiction is obtained.  \QED

This theorem can leads us to the following conclusion: if the complete matching exists for any 
possible removal of  $q-1$ vertexes from Koenig representation of graph $\langle \Gamma (W
\cup \{d_i\}), W \cup \{d_i\} \rangle$ where $W \cup \{ d_i \} \subseteq D$ this graph can 
be taken as the subgraph for skeleton developed with the use of greedy algorithm. 
For example removal of hypergraph   $\tilde \Gamma$~(see Figure~1 and Figure~2) vertexes
\{3,5\}, or \{3,6\}, etc. from Koening representation leads to complete matching 
absence, i.e. edges set \{c,d,e,f\} of hypergraph $\tilde \Gamma$ is independent. If suppose that the
hypergraph $\tilde \Gamma$ is the current subgraph of skeleton developed with the use of 
greedy algorithm then the previously added  edge is to be thrown off (e.g. the edge $f$ can be taken
for that purpose).

{\bf Theorem 2.} {\it Optimal hypergraph $\Gamma = \langle Z, D
\rangle$ skeleton development algorithm is the polynomial one.}

{\bf Proof.} Theorem~1 provides us the following conclusion: to
find independence of edges subset  $A \subseteq D$ of hypergraph
 $\langle \Gamma A, A \rangle$ it is necessary to perform complete matching existence
 check of Koenig representation of hypergraph $\langle \Gamma A, A \rangle$ for $k =
{|\Gamma A| \choose q-1} $ times. The complete matching
development algorithm is polynomial one~[2] and the total amount
of complete matching algorithm call in greedy algorithm is lower
than $|D|^q$, therefore theorem statement is true. \QED

In case of direct implementation of greedy algorithm to develop
the hypergraph skeleton the number of used operations can be
considerably reduced; it can be concluded from the following
statement.

{\bf Theorem 3.} {\it Let  hypergraph  $ \langle Z', D'\rangle$ of
hypergraph $ \langle Z, D \rangle$ is hyperforest of ($ Z' \subset
Z, D'\subset D$) и $a \in D \setminus D'$.

Then

1) If  $|\Gamma a \cap \Gamma D'| < q$, then in hypergraph
$\langle Z'+ \Gamma a, D'+ a \rangle$ edges set  $(D'+ a)$ is
independent one and therefore it is unnecessary to check the
complete matching existence in hypergraph Koenig representation;

2) If $|\Gamma a \cap \Gamma D'| = q$ then to check edges set
$(D'~+~a)$ independence \ it is enough \ to \ check \ complete \
matching \ existence \ in \ hypergraph \ Koenig representation
$\langle \Gamma (D'+a), D' + a \rangle$ while removing any
 $q-1$ vertexes of the subset~$\Gamma a$.}


{\bf Proof.} Genuinely, in the first case the edges subset $D'$ is
independent by definition, i.e. for any nonvacuous  $A \subset D'$
the inequality  $\left| \Gamma A \right| \ge \left| A \right| + q
- 1$ is true. When edge  $a$ is added to subset $D'$ the vertex
amount is increased at least for one i.e. $|\Gamma (D' + a)| \ge
|\Gamma D'| + 1$. Adding the left and right parts of these
inequalities to one another we can obtain that for any nonvacuous
 $A \subset D'+ a $ the formula $\left|
\Gamma A + a \right| \ge \left| A + a \right| + q - 1$ is true.
Therefore first part of theorem 3 is proved.

Let us consider two cases to prove the second part of the theorem
3.

1. Let us suppose that having the edge  $a$ added to the set $D'$
the total amount of combined vertexes is equal to $q$, i.e.
$|\Gamma a \cap \Gamma D'| = q$ and obtained hypergraph
 $\langle \Gamma (D' + a), D' + a
\rangle$ remained the one of hyperforest type. Therefore,
according to the theorem~1, while removing any $q-1$ vertexes of
vertex set $\langle \Gamma (D' + a) \rangle$ of hypergraph Koenig
representation the complete matching is exist in this
representation.  Thus, matching exist in case while any $q-1$
vertexes of subset $\Gamma a$ are chosen as removable ones.

2. Let us consider that having the edge  $a$ added to the set $D'$
the $|\Gamma a \cap \Gamma D'| = q$ as it was in previous case;
however dependent edges subset  $D'' \subseteq D'$ is developed in
obtained hypergraph  $\langle \Gamma (D' + a), D' + a \rangle$.

Let $D' = T_1 + T_2 + \ldots + T_b$ is set $D'$ partition
corresponds the  skeleton connection component of hypergraph
 $ \langle Z', D'\rangle$.
Because an edge  $a \in D''$  has in this case the unique  $i \in
1:b$ such as $D'' \subseteq T_i + a$ then $a$ is a link of the
corresponding hypergraph skeleton connection component  $T_i$.
Condition $|\Gamma D''| - (q - 1) < |D''|$ is true by definition
for $D''$ subset. Therefore it can be concluded that removing any
$q-1$ vertexes of $D''$ the Hall theorem conditions of subgraph
$\langle \Gamma D'', D'' \rangle$ Koenig representation are
failed.
 Thus, while checking the independence any $q-1$ vertexes belonged to $\Gamma a$ can be taken as removable ones. 
 \QED

While removing the $q-1$ vertexes of subset  $\Gamma a$ of
hypergraph Koenig space $\langle \Gamma (D' + a), D' + a \rangle$
the complete matching is absent then it can be concluded that this
hypergraph includes the dependent set.  Otherwise obtained graph
is the hyperforest and the next edge of $\langle \Gamma D, D
\rangle$ can be chosen and analyzed then.

It is obvious that utilizing this theorem results at the third
step of the greedy algorithm the possible amount of developed
skeleton Koenig representation complete matching separation
algorithm call can be decreased.

\vspace{3mm} 
\centerline{REFERENCES} 
\vspace{1mm}

1. Aigner M. Theory of Combinations. --- M., Mir, 1982. --- 556~p

2. Lipsky V. Combinatorics for Programmers.  --- M., Mir, 1988. --- 213~p.

3. Sushkov Yu.A. (1,q)-combinations. --- Spb., Vestnik LGU. 1975,
N~4.--- P.~50-55.

4. Sushkov Yu.A. Hypergraphs Connection. --- Computation Technics
and Cybernetics Problems. SPb., LGU. 1984, issue~20. --- P.~87--96.

5. Sushkov Yu.A. Matroids and Hypergraphs. --- Vestnik LGU.
Mathematics.  Mechanics. Astronomy.  SPb., LGU, 1985, N 22. --- P.~42--46.

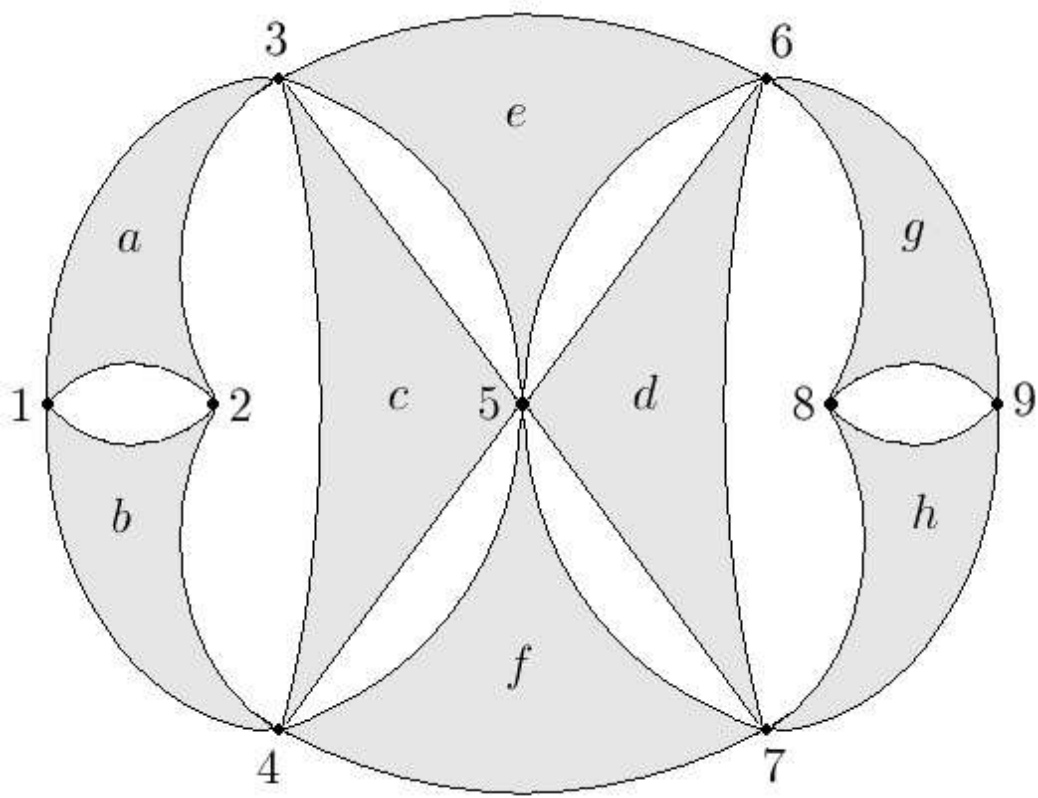
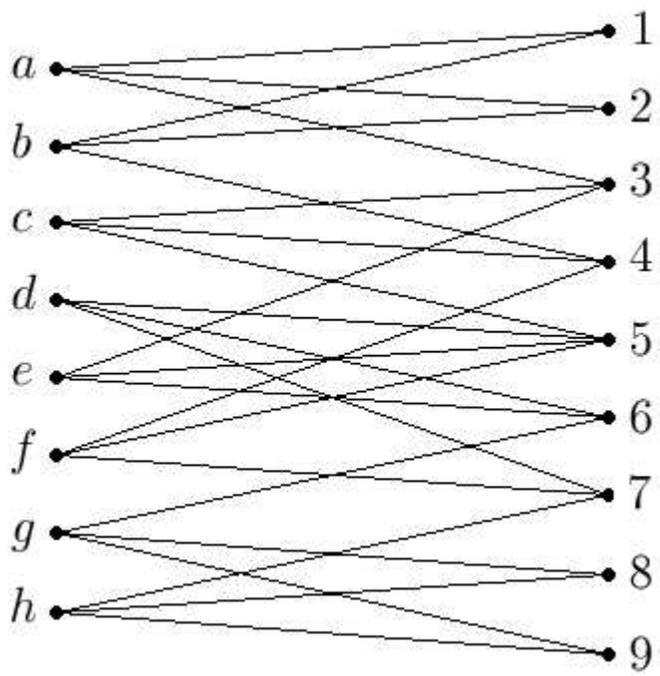
\end{document}